\documentclass[11pt]{article}
\usepackage[francais]{babel}
\usepackage[latin1]{inputenc}
\usepackage[T1]{fontenc}
\usepackage{a4,amsmath,amsfonts}
\title{Variétés abéliennes et théorème de Minkowski-Hlawka}
\author{Pascal Autissier}
\begin{document}

\maketitle

\newcommand{\D}{\displaystyle}

{\bf Abstract:} A classical theorem of Minkowski and Hlawka states that there
exists a lattice in $\mathbb{R}^n$ with packing density at least $2^{1-n}$.
Buser and Sarnak proved the analogue of this result in the context of complex
abelian varieties. Here we give an improvement of this analogue; this shows a
conjecture of Muetzel.\\

{\bf Résumé:} Un théorème classique de Minkowski et Hlawka montre l'existence
d'un réseau de $\mathbb{R}^n$ à densité d'empilement $\ge2^{1-n}$. Buser et
Sarnak ont établi l'analogue de ce résultat dans le cadre des variétés
abéliennes complexes. On donne ici une amélioration de cet analogue; cela
prouve une conjecture de Muetzel.\\

{\it 2010 Mathematics Subject Classification:} 11H31, 14K20.\\

\section{Introduction}

Un problème important de la Géométrie des Nombres est l'étude de la plus grande
densité d'empilement des réseaux d'un espace euclidien $E$ de dimension
$n\ge1$. Rappelons que si $\Gamma$ est un réseau de $E$ de premier minimum
$\lambda_1(\Gamma)$, sa densité d'empilement $\Delta(\Gamma)$ est la densité de
l'empilement de boules de rayon $\D\frac{\lambda_1(\Gamma)}{2}$ centrées en les
points de $\Gamma$. Le problème en question est alors d'estimer
$\Delta_n=\sup\{\Delta(\Gamma)\ ;\ \Gamma\mbox{ réseau de }E\}$.\\

La valeur exacte de $\Delta_n$ n'est établie que pour certains petits entiers
$n$, et le comportement asymptotique de $\Delta_n$ lorsque $n$ tend
vers $+\infty$ reste mystérieux. Le fameux théorème de Minkowski-Hlawka donne
l'inégalité $\D\Delta_n\ge\frac{1}{2^{n-1}}$ pour tout $n\ge1$. La meilleure
minoration générale connue (à la valeur de la constante $c$ près) est due à
Rogers \cite{Roge} en 1947:\\

{\bf Théorème (Rogers):} {\it En posant $\D c=\frac{2}{e}$, on a
$\D\Delta_n\ge\frac{cn}{2^n}$ pour tout entier $n\ge1$.}\\

La constante $c$ a été améliorée successivement par Davenport-Rogers
\cite{DaRo}, Ball \cite{Ball}, Venkatesh \cite{Venk}. Ce dernier obtient
$\D\Delta_n\ge\frac{c'n}{2^n}$ pour tout $n$ assez grand, avec $c'=65963$. Dans
le même article, il montre également le résultat suivant:\\

{\bf Théorème (Venkatesh):} {\it Il existe une infinité d'entiers $n\ge3$
vérifiant l'inégalité $\D\Delta_n\ge\frac{n}{2^{n+1}}\ln\ln n$.}\\

\`A titre de comparaison, la meilleure majoration générale connue est celle de
Kabatyanskii et Levenshtein \cite{KaLe}, qui est de la forme $\Delta_n\le C^n$
pour tout $n$ assez grand, avec $C=0,661$.\\

On s'intéresse dans ce travail à l'analogue de ce problème dans le contexte des
variétés abéliennes complexes.

Soit $(A;L)$ une variété abélienne complexe de dimension $g\ge1$,
principalement polarisée. Posons $T_A=\Gamma(A;\Omega_{A/\mathbb{C}})^\vee$ et
désignons par $\Gamma_A$ le réseau des périodes de $A$ (on a donc un
isomorphisme $A(\mathbb{C})\simeq T_A/\Gamma_A$ de groupes analytiques). La
polarisation $L$ induit une forme de Riemann $\bigl<\ ;\ \bigr>$ sur $T_A$,
{\it i.e.} un produit scalaire hermitien sur $T_A$ tel que
${\rm Im}\bigl<\gamma_1;\gamma_2\bigr>\in\mathbb{Z}$ pour tout
$(\gamma_1;\gamma_2)\in\Gamma_A^2$.\\

{\bf Définition:} Le {\bf volume d'injectivité} $V(A;L)$ de $(A;L)$ est la
densité d'empilement du réseau $\Gamma_A$ dans l'espace euclidien
$(T_A;{\rm Re}\bigl<\ ;\ \bigr>)$. De manière équivalente, $V(A;L)$ est le plus
grand volume d'une boule ouverte de $T_A$ qui s'injecte dans $A$ (le volume est
normalisé de sorte que $\Gamma_A$ soit de covolume 1).\\

\'Etant donné un entier $g\ge1$, on note ${\cal A}_g$ l'espace de modules
(grossier) des schémas abéliens de dimension relative $g$ et principalement
polarisés. On cherche ici à minorer
$V_g=\sup\{V(A;L)\ ;\ (A;L)\in{\cal A}_g(\mathbb{C})\}$. On a trivialement
l'inégalité $V_g\le\Delta_{2g}$. Buser et Sarnak \cite{BuSa}
ont montré l'équivalent du théorème de Minkowski-Hlawka dans ce cadre:\\

{\bf Théorème (Buser, Sarnak):} {\it Pour tout entier $g\ge1$, on a
$\D V_g\ge\frac{1}{2^{2g-1}}$.}\\

On se propose ici de raffiner ce résultat pour certaines valeurs de $g$. On
désigne par $\varphi$ la fonction indicatrice d'Euler.\\

{\bf Théorème 1.1:} {\it Soit $m$ un entier $\ge3$; posons $g=\varphi(m)$. On a
la minoration $\D V_g\ge\frac{m}{4^g}$.}\\

Tirons quelques conséquences du théorème 1.1. Muetzel conjecture dans
\cite{Muet} (voir sa conjecture 2.5) l'estimation $4^gV_g\ge2g$ lorsque $g$ est
une puissance de 2. On obtient ici une version améliorée de cet énoncé:\\

{\bf Corollaire 1.2:} {\it $(\alpha)$ Soit $g$ un entier $\ge2$ qui est une
puissance de 2. On a $4^gV_g\ge3g$.

$(\beta)$ Soit $g$ un entier de la forme $g=\varphi(n)$ pour un entier $n\ge3$.
On a $4^gV_g\ge2g+2$.}\\

Ce résultat suggère que $V_g$ admet peut-être une minoration semblable à celle
de Rogers:\\

{\bf Question:} Si $g$ est un entier $\ge2$, a-t-on $4^gV_g\ge2g+2$ ?\\

On déduit aussi de l'énoncé 1.1 un analogue (ou plutôt un raffinement, au vu de
la majoration $V_g\le\Delta_{2g}$) du théorème de Venkatesh:\\

{\bf Corollaire 1.3:} {\it Notons ici $\gamma$ la constante d'Euler. Il existe
un réel $c$ et une infinité d'entiers $g$ tels que
$4^gV_g\ge e^\gamma g\ln\ln g-cg$.}\\

La démonstration du théorème 1.1 repose sur un argument de valeur moyenne sur
un certain espace de réseaux, inspiré de la méthode classique de
Minkowski-Hlawka. Plus précisément, on considère des réseaux
$\Gamma\subset\mathbb{C}^g$ munis d'une action libre d'un groupe cyclique
d'ordre $m$ et tels que $\mathbb{C}^g/\Gamma$ soit une variété abélienne ayant
un anneau d'endomorphismes de rang $\D\ge\frac{g}{2}$ (voir remarque 4.1).\\

Je remercie Renaud Coulangeon pour m'avoir fourni la référence \cite{Venk}. Je
remercie également Gaël Rémond pour ses commentaires sur ce travail.\\

\section{Rappels}

Soit $n$ un entier $\ge1$. Munissons $\mathbb{R}^n$ de la norme euclidienne
usuelle, notée $\|\ \|$, et de la mesure de Lebesgue. On désigne par $v_n$ le
volume de la boule unité de $\mathbb{R}^n$, de sorte que
$\D v_n=\frac{\pi^{n/2}}{(n/2)!}$.

Soit $\Gamma$ un réseau de $\mathbb{R}^n$ de covolume 1. Le premier minimum de
$\Gamma$ est par définition le réel
$\D\lambda_1(\Gamma)=\min_{\gamma\in\Gamma-\{0\}}\|\gamma\|$. La densité d'empilement
$\Delta(\Gamma)$ de $\Gamma$ est la densité de l'empilement de boules de rayon
$\D\frac{\lambda_1(\Gamma)}{2}$ centrées en les points de $\Gamma$. On a donc
la formule $\D\Delta(\Gamma)=\frac{v_n}{2^n}\lambda_1(\Gamma)^n$.\\

Soit $g$ un entier $\ge1$. On considère maintenant $\mathbb{C}^g$ muni du
produit scalaire hermitien standard ($\mathbb{C}$-linéaire à gauche par
convention) noté $\bigl<\ ;\ \bigr>$, et de la mesure de Lebesgue.

Soit $\Gamma$ un réseau (de rang $2g$) de $\mathbb{C}^g$ tel que
${\rm Im}\bigl<\gamma;\gamma'\bigr>\in\mathbb{Z}$ pour tout
$(\gamma;\gamma')\in\Gamma^2$. Le tore $\mathbb{C}^g/\Gamma$ est alors une
variété abélienne complexe, et la forme de Riemann $\bigl<\ ;\ \bigr>$ induit
une polarisation $L$ sur $\mathbb{C}^g/\Gamma$. Posons
$\Gamma'=\{\gamma'\in\mathbb{C}^g\ |\ \forall\gamma\in\Gamma\ {\rm Im}\bigl<\gamma;\gamma'\bigr>\in\mathbb{Z}\}$; c'est un réseau de $\mathbb{C}^g$ contenant
$\Gamma$.  Si $\Gamma'=\Gamma$, alors la polarisation $L$ est principale et
$\Gamma$ est de covolume 1.

Pour des précisions sur les polarisations, on pourra consulter \cite{BiLa}
pages 69-74.\\

\section{Préliminaires}

Soient $m$ un entier $\ge3$ et $\zeta\in\overline{\mathbb{Q}}$ une racine
$m$-ième primitive de l'unité. On pose $g=\varphi(m)$ et on plonge
$K=\mathbb{Q}(\zeta)$ dans $\mathbb{C}^g$ via les $g$ plongements complexes de
$K$. On désigne par $E$ le sous-espace vectoriel \underline{réel} de
$\mathbb{C}^g$ engendré par $K$. Observons que $E$ est en fait un sous-anneau
de $\mathbb{C}^g$ identifié à $K\otimes_\mathbb{Q}\mathbb{R}$, et que
$\mathbb{C}^g=E\oplus iE$.\\

Pour tout $z=(z_1;\cdots;z_g)\in\mathbb{C}^g$, on note
$\bar{z}=(\bar{z_1};\cdots;\bar{z_g})$ son conjugué complexe. Remarquons que
$K$, $O_K$ et $E$ sont stables par cette conjugaison (car
$\bar{\zeta}=\zeta^{-1}$). En outre, le produit scalaire hermitien standard
$\bigl<\ ;\ \bigr>$ sur $\mathbb{C}^g$ est à valeurs réelles sur $E^2$, donc
fait de $E$ un espace euclidien. Plus précisément, on a
$\bigl<a;b\bigr>={\rm Tr}(a\bar{b})$ pour tout $(a;b)\in K^2$, où
${\rm Tr}:K\rightarrow\mathbb{Q}$ désigne la trace.\\

Notons $G$ le sous-groupe (cyclique d'ordre $m$) de $O_K^*$ engendré par
$\zeta$. On définit une action $\mathbb{R}$-linéaire de $G$ sur
$\mathbb{C}^g=E\oplus iE$ par $g*(x+iy)=gx+i\bar{g}y$ pour tout
$(g;x;y)\in G\times E^2$. Cette action est libre sur $\mathbb{C}^g-\{0\}$.\\

{\bf Lemme 3.1:} {\it Soient $g\in G$ et $z\in\mathbb{C}^g$. On a
$\|g*z\|=\|z\|$.}\\

{\it Démonstration:} \'Ecrivons $z=x+iy$ avec $x$ et $y$ dans $E$. On a alors
$\|z\|^2=\|x\|^2+2{\rm Re}(i\bigl<y;x\bigr>)+\|y\|^2=\|x\|^2+\|y\|^2$. En
utilisant l'égalité $|\sigma(g)|=1$ pour tout plongement
$\sigma:K\hookrightarrow\mathbb{C}$, on obtient $\|gx\|^2=\|x\|^2$. On en
conclut que

$\|g*z\|^2=\|gx+i\bar{g}y\|^2=\|gx\|^2+\|\bar{g}y\|^2=\|x\|^2+\|y\|^2=\|z\|^2$.
$\square$\\

Désignons par $I$ la codifférente de $O_K$, c'est-à-dire l'idéal fractionnaire
défini par $I=\{a\in K\ |\ {\rm Tr}(aO_K)\subset\mathbb{Z}\}$. On sait que
$\{b\in K\ |\ {\rm Tr}(bI)\subset\mathbb{Z}\}=O_K$ et que $I$ est un réseau de
$E$. En outre, on voit aisément que $I$ est stable par conjugaison complexe.\\

{\bf Lemme 3.2:} {\it Soit $f:E\rightarrow E$ une application
$\mathbb{R}$-linéaire symétrique. Soit $r$ un réel $>0$. Considérons le réseau
$$\Gamma=rI+\Bigl(rf+\frac{i}{r}{\rm Id}\Bigr)(O_K)=\Bigl\{ra+rf(b)+\frac{i}{r}b\ ;\ (a;b)\in I\times O_K\Bigr\}\mbox{ de }\mathbb{C}^g.$$

$(\alpha)$ Le tore $\mathbb{C}^g/\Gamma$ est une $\mathbb{C}$-variété abélienne
principalement polarisée.

$(\beta)$ Soit $x\in E$. Si $f$ est l'application $\begin{array}{c}
E\rightarrow E\\
y\mapsto x\bar{y}\\
\end{array}$, alors $f$ est symétrique, et $\Gamma$ est stable sous l'action de
$G$.}\\

{\it Démonstration:} $(\alpha)$ Vérifions que $\bigl<\ ;\ \bigr>$ est une forme
de Riemann. Soient $(a;a';b;b')\in I^2\times O_K^2$; posons
$\D\gamma=ra+rf(b)+\frac{i}{r}b$ et $\D\gamma'=ra'+rf(b')+\frac{i}{r}b'$. On a
alors
$${\rm Im}\bigl<\gamma;\gamma'\bigr>=\bigl<b;a'\bigr>+\bigl<b;f(b')\bigr>-\bigl<a;b'\bigr>-\bigl<f(b);b'\bigr>={\rm Tr}(b\bar{a'})-{\rm Tr}(a\bar{b'})\in\mathbb{Z}\ .\qquad(*)$$

Le quotient $\mathbb{C}^g/\Gamma$ est donc une variété abélienne polarisée.\\

Montrons que la polarisation est principale. Posons
$\Gamma'=\{\gamma'\in\mathbb{C}^g\ |\ \forall\gamma\in\Gamma\ {\rm Im}\bigl<\gamma;\gamma'\bigr>\in\mathbb{Z}\}$. Le réseau $\Gamma$ est d'indice fini dans
$\Gamma'$, ce qui implique $\Gamma'\subset\Gamma\otimes_\mathbb{Z}\mathbb{Q}$.
Soit $\gamma'\in\Gamma'$; il existe donc $(a';b')\in K^2$ tel que
$\D\gamma'=ra'+rf(b')+\frac{i}{r}b'$.

Pour tout $b\in O_K$, on a $\D{\rm Tr}(a'b)={\rm Im}\bigl<rf(\bar{b})+\frac{i}{r}\bar{b};\gamma'\bigr>\in\mathbb{Z}$ par un calcul similaire à $(*)$, donc
$a'\in I$. De même, on a
${\rm Tr}(b'a)=-{\rm Im}\bigl<r\bar{a};\gamma'\bigr>\in\mathbb{Z}$ pour tout
$a\in I$, ce qui donne $b'\in O_K$. D'où $\gamma'\in\Gamma$. On a bien prouvé
que $\Gamma'=\Gamma$, {\it i.e.} que la polarisation est principale.\\

$(\beta)$ On suppose maintenant que $f$ est l'application qui à $y\in E$
associe $x\bar{y}$. Soit $(y;z)\in E^2$. En écrivant par coordonnées les
vecteurs $x$, $y$ et $z$ de $\mathbb{C}^g$, on trouve
$$\bigl<y;f(z)\bigr>=\bigl<f(z);y\bigr>=\sum_{k=1}^gx_k\bar{z_k}\bar{y_k}=\bigl<f(y);z\bigr>\quad.$$
Ainsi $f$ est-elle symétrique. Soient $g\in G$ et $(a;b)\in I\times O_K$; on
pose $\D\gamma=ra+rf(b)+\frac{i}{r}b$. On a
$\D g*\gamma=g(ra+rx\bar{b})+\frac{i}{r}\bar{g}b=rga+rf(\bar{g}b)+\frac{i}{r}\bar{g}b\in\Gamma$. D'où la $G$-stabilité de $\Gamma$. $\square$\\

On choisit un domaine fondamental $F$ du réseau $I$ dans $E$, et un domaine
fondamental $F'$ de $O_K$. Désignons par $\nu$ la mesure de Lebesgue sur $E$,
normalisée de sorte que la boule unité dans $E$ soit de volume $v_g$. Avec
cette normalisation, la mesure image de $\nu\times\nu$ par l'application
$h:\begin{array}{c}
E^2\rightarrow\mathbb{C}^g\\
(x;y)\mapsto x+iy\\
\end{array}$ s'identifie à la mesure de Lebesgue usuelle sur $\mathbb{C}^g$.\\

Regardons le réseau $\Gamma$ du lemme 3.2 avec $f=0$: ce réseau
$\Gamma=h(rI\times\frac{1}{r}O_K)$ est de covolume $\nu(rF)\nu(\frac{1}{r}F')$
d'une part, et de covolume 1 d'autre part. On en déduit que
$\nu(F)\nu(F')=1$.\\

{\bf Lemme 3.3:} {\it Soient $b\in O_K-\{0\}$, $z\in\mathbb{C}^g$ et
$r\in\mathbb{R}^*_+$. Soit $\chi:\mathbb{C}^g\rightarrow\mathbb{C}$ une
application telle que la fonction $E\rightarrow\mathbb{C}$ qui à $x$ associe
$\chi(x+z)$ soit intégrable et à support compact. On a l'égalité
$$\int_F\sum_{a\in I}\chi(ra+rxb+z){\rm d}\nu(x)=\frac{1}{r^g}\int_E\chi(x+z){\rm d}\nu(x)\quad.$$}

{\it Démonstration:} L'application $\psi:E\rightarrow\mathbb{C}$ qui à $x$
associe $\D\sum_{a\in I}\chi(ra+rx+z)$ est $I$-périodique, et la multiplication
par $b$ induit un revêtement $E/I\rightarrow E/I$. La formule de changement
(linéaire) de variable donne
$\D\int_{E/I}\psi(bx){\rm d}\nu(x)=\int_{E/I}\psi(x){\rm d}\nu(x)$. On en déduit
que
$$\int_F\psi(xb){\rm d}\nu(x)=\sum_{a\in I}\int_F\chi(ra+rx+z){\rm d}\nu(x)=\int_E\chi(rx+z){\rm d}\nu(x)=\frac{1}{r^g}\int_E\chi(x+z){\rm d}\nu(x)\ .$$

D'où le résultat. $\square$\\

\section{Démonstration du théorème 1.1}

On conserve les notations de la partie 3. Soit $\varepsilon$ un réel vérifiant
$0<\varepsilon<m$. On va prouver l'existence d'une $\mathbb{C}$-variété
abélienne principalement polarisée $(A;L)$ telle que
$\D V(A;L)>\frac{m-\varepsilon}{4^g}$. Lorsque $z\in\mathbb{C}^g$, on pose
$\chi(z)=1$ si $v_{2g}\|z\|^{2g}\le m-\varepsilon$ et $\chi(z)=0$ sinon. Pour
tout réel $r>0$, on pose
$$J(r)=\frac{\nu(F')}{r^g}\sum_{b\in O_K-\{0\}}\int_E\chi\Bigl(x+\frac{i}{r}b\Bigr){\rm d}\nu(x)\quad.$$

$J(r)$ est une somme de Riemann, donc $J(r)$ converge vers $\D\int_E\int_E\chi(x+iy){\rm d}\nu(x){\rm d}\nu(y)=m-\varepsilon$ lorsque $r$
tend vers $+\infty$. En particulier, il existe un réel $r_0>0$ vérifiant
$J(r_0)<m$ et $v_{2g}(r_0\lambda_1(I))^{2g}>m-\varepsilon$.\\

On définit une application $N:\mathbb{R}^g\rightarrow\mathbb{N}$ par la formule
$$\forall x\in\mathbb{R}^g\qquad N(x)=\sum_{b\in O_K-\{0\}}\sum_{a\in I}\chi\Bigl(r_0a+r_0x\bar{b}+\frac{i}{r_0}b\Bigr)\quad.$$

Grâce au lemme 3.3, on a la relation
$\D\int_FN(x){\rm d}\nu(x)=\frac{J(r_0)}{\nu(F')}$. Il existe donc un $x_0\in F$
tel que $\D N(x_0)\le\frac{J(r_0)}{\nu(F)\nu(F')}<m$.\\

Notons $f:E\rightarrow E$ l'application $y\mapsto x_0\bar{y}$ et posons
$\D\Gamma=r_0I+\Bigl(r_0f+\frac{i}{r_0}{\rm Id}\Bigr)(O_K)$. D'après le lemme
3.2, le quotient $A=\mathbb{C}^g/\Gamma$ est une variété abélienne
naturellement munie d'une polarisation $L$ principale.\\

Il reste à minorer $V(A;L)$. Soit
$\D\gamma=r_0a+r_0f(b)+\frac{i}{r_0}b\in\Gamma-\{0\}$. Si $b\neq0$, le lemme 3.1
permet d'obtenir $\D m\chi(\gamma)=\sum_{g\in G}\chi(g*\gamma)\le N(x_0)<m$, ce
qui implique $v_{2g}\|\gamma\|^{2g}>m-\varepsilon$. Et si $b=0$, on a directement
$v_{2g}\|\gamma\|^{2g}\ge v_{2g}(r_0\lambda_1(I))^{2g}>m-\varepsilon$.

On a donc montré que
$\D V(A;L)=\frac{v_{2g}}{4^g}\lambda_1(\Gamma)^{2g}>\frac{m-\varepsilon}{4^g}$.
D'où le résultat. $\square$\\

{\bf Remarque 4.1:} Considérons le sous-corps totalement réel
$K'=\mathbb{Q}(\zeta+\zeta^{-1})$ de $K$. Soit $b\in O_{K'}$; notons
$j(b):\mathbb{C}^g\rightarrow\mathbb{C}^g$ la multiplication par $b$. Le réseau
$\Gamma$ du lemme 3.2.$\beta$ est alors stable par $j(b)$. Cela induit un
endomorphisme $j(b):\mathbb{C}^g/\Gamma\rightarrow\mathbb{C}^g/\Gamma$. On
vient ainsi de construire un morphisme injectif d'anneaux
$j:O_{K'}\hookrightarrow{\rm End}(\mathbb{C}^g/\Gamma)$.

En posant $W_{K'}=\sup\{V(A;L)\ ;\ (A;L)\in{\cal A}_g(\mathbb{C})\mbox{ et }O_{K'}\mbox{ s'injecte dans }{\rm End}(A)\}$, la démonstration ci-dessus prouve en
fait l'inégalité $\D W_{K'}\ge\frac{m}{4^g}$.\\

\section{Démonstration des corollaires}

{\it Démonstration de 1.2:} $(\alpha)$ Il suffit d'appliquer le théorème 1.1
avec $m=3g$.\\

$(\beta)$ Si $n$ est pair sans être une puissance de 2, on a
$4^gV_g\ge n\ge2g+2$. Si $n$ est impair, on a
$\varphi(2n)=g$ donc $4^gV_g\ge2n\ge2g+2$. $\square$\\

{\it Démonstration de 1.3:} Soit $x$ un réel $\ge3$. Prenons $m$ égal au
produit des nombres premiers $p\le x$, et posons $g=\varphi(m)$. D'après le
théorème de Mertens (voir \cite{Tene} page 17), on a
$\D\frac{m}{g}=\prod_{p\le x}\Bigl(1-\frac{1}{p}\Bigr)^{-1}=e^\gamma\ln x+O(1)$.
Les estimations de Tchébychev donnent $\ln g\le\ln m=O(x)$, ce qui implique
$\ln\ln g\le\ln x+O(1)$. On en conclut que
$4^gV_g\ge m\ge e^\gamma g\ln\ln g+O(g)$. $\square$\\

\ \\

{\small Pascal Autissier. I.M.B., université de Bordeaux, 351, cours de la
Libération, 33405 Talence cedex, France.

pascal.autissier@math.u-bordeaux1.fr}

\end{document}